\documentclass[leqno,12pt]{amsart}
\usepackage{latexsym}
\usepackage{mathrsfs}
\usepackage{amsmath}
\usepackage{amssymb}
\usepackage[bookmarks]{hyperref}
\usepackage{pstricks,pst-plot}

\newtheorem{theorem}{Theorem}[section]
\newtheorem{lemma}[theorem]{Lemma}
\newtheorem{corollary}[theorem]{Corollary}
\newtheorem{proposition}[theorem]{Proposition}

\theoremstyle{definition}

\newtheorem{definition}[theorem]{Definition}

\newcommand{\C}{\mathbb{C}}

\newcommand{\R}{\mathbb{R}}

\newcommand{\U}{\mathbb{U}}

\newcommand{\bthm}{\begin{theorem}}
\newcommand{\ethm}{\end{theorem}}
\newcommand{\blem}{\begin{lemma}}
\newcommand{\elem}{\end{lemma}}
\newcommand{\bcor}{\begin{corollary}}
\newcommand{\ecor}{\end{corollary}}
\newcommand{\bprop}{\begin{proposition}}
\newcommand{\eprop}{\end{proposition}}
\newcommand{\bdefn}{\begin{definition}}
\newcommand{\edefn}{\end{definition}}
\newcommand{\bpf}{\begin{proof}}
\newcommand{\epf}{\end{proof}}

\def\vep {\varepsilon}
\def \sm {\setminus}
\def \ob {\overline B}

\def\h#1{\widehat {#1}}

\def \ab#1{{[a_{#1}, b_{#1}]}}
\def \openab#1{(a_{#1},b_{#1})}

\def \Int {{\rm Int}}
\def\I {[0,1]}

\def\albeta {(\alpha_j,\beta_j)}
\def\tsigma {\tilde\sigma}
\def\diam {{\rm diam}}

\def\itemskip {\vskip 4pt plus 2 pt minus 2 pt}

\begin{document}
\title{Polynomial Hulls of Arcs and Curves II}
\author{Alexander J. Izzo}
\address{Department of Mathematics and Statistics, Bowling Green State University, Bowling Green, OH 43403}
\email{aizzo@bgsu.edu}
\thanks{The author was supported by NSF Grant DMS-1856010.}


\begin{abstract}
We prove that if a compact set $E$ in $\C^N$ is contained in an arc $J$, then there is a choice of $J$ whose polynomial hull $\h J$ is $J\cup \h E$.  This strengthens an earlier result of the author.  We also correct an inaccuracy in the statement, and fill a gap in the proof, of that earlier result. 
\end{abstract}

\maketitle

%
%
%
%

\section{Introduction}

The purpose of this paper is to strengthen results in the author's earlier paper \cite{I2021} and to address a gap in that paper.  Our main result shows that every compact set $E$ that is contained in an arc in $\C^N$ is contained in an arc whose polynomial hull is no larger than it is obviously forced to be by virtue of containing the set $E$.

\bthm\label{general-theorem}
Let $E$ be a compact set that is contained in an arc in $\C^N$.
Then some arc $J$ in $\C^N$ that contains $E$ is such that $\h J=J \cup \h E$.  Furthermore, $J$ can be chosen to lie in an arbitrary connected neighborhood of $E$ and such that each component of $J\sm E$ is a $C^\infty$-smooth open arc.  In addition, $J$ can be taken to be a 
simple closed curve rather than an arc provided $N\geq 2$.
\ethm


The special case when the set $E$ is polynomially convex yields the following.

\bcor\label{corollary}
Each compact polynomially convex set $E$ that is contained in an arc in $\C^N$ is contained in a polynomially convex arc $J$ that can be chosen to lie in an arbitrary connected neighborhood of $E$.  Furthermore, $J$ can be chosen so that each component of $J\sm E$ is a $C^\infty$-smooth open arc.
With that choice, if $P(E)=C(E)$, then $P(J)=C(J)$.
In addition, $J$ can be taken to be a 
simple closed curve rather than an arc provided $N\geq 2$.
\ecor

We recall here some standard 
terminology and notation already used above.
By definition an \emph{arc} is a space homeomorphic to the closed unit interval, and a \emph{simple closed curve} is a space homeomorphic to the unit circle.  For convenience  we will also use the term \emph{arc} (or \emph{simple closed curve}) to refer to a topological embedding whose domain is the closed unit interval (or the unit circle).
A mapping that is referred to as a $C^\infty$-\emph{smooth arc} will be required to be an immersion, i.e., to have nowhere vanishing derivative.
Throughout the paper, $m$ will denote $2N$-dimensional Lebesgue measure on $\C^N$.
For a compact set $X$ in $\C^N$, we denote by $C(X)$ the space of all continuous complex-valued functions on $X$ and by
$P(X)$ the uniform closure in $C(X)$ of the polynomials in the complex coordinate functions $z_1,\ldots, z_N$.
The \emph{polynomial hull} 
$\h X$ of $X$ is defined by
$$\h X=\{z\in\C^N:|p(z)|\leq \max_{x\in X}|p(x)|\ 
\mbox{\rm{for\ all\ polynomials}}\ p\}.$$
The set $X$ is said to be \emph{polynomially convex} if $\h X = X$.  

For the special case in which the set $E$ is totally disconnected, a result along the lines of Theorem~\ref{general-theorem} was presented in \cite[Theorem~1.3]{I2021} and used there to establish the existence of arcs with certain properties.  However, the proof given there seems to have a gap.  It was asserted that given the existence of an arc that contains the totally disconnected, compact set $E$, \lq\lq one can show that there is such an arc $J_0$ with the additional
property that the closure of each component of $J_0\sm E$ is a $C^\infty$-smooth arc\rq\rq.  We will prove in the present paper that there is an arc through $E$ such that each component of $J_0\sm E$ is a $C^\infty$-smooth open arc, but the author does not know
whether these open arcs can be chosen so that their closures are smooth (closed) arcs.  That stronger condition was used in the proof given for \cite[Theorem~1.3]{I2021} in that it gave that $J_0\sm E$ was contained in a countable union of \emph{compact} one-dimensional smooth manifolds-with-boundary and thus made it possible to use the stability of smooth embeddings in the $C^1$-topology.  

The (flawed) approach used in \cite{I2021} can be adapted to give a correct proof using the stability of smooth embeddings in the strong topology (also known as the fine topology or the Whitney topology).  (See \cite[p.~35]{Hirsch:1976} for the definition.)  This, however, makes the details somewhat more complicated.  We will instead use a different approach.  Theorem~\ref{general-theorem} is closely related to \cite[Theorem~1.1]{IS} and can, in fact, be regarded as a generalization of that theorem.  The proof we will give for Theorem~\ref{general-theorem} will use results from \cite{IS}.

Fortunately the gap in the proof of \cite[Theorem~1.3]{I2021} has very little effect on the applications in \cite{I2021}.  The proof of \cite[Theorem~1.1]{I2021}, which gives the existence of arcs and simple closed curves in $\C^3$ having \lq\lq hull with dense invertibles\rq\rq\ goes through unchanged except for invoking Theorem~\ref{general-theorem} of the present paper in place of \cite[Theorem~1.3]{I2021}.

Likewise \cite[Theorem~1.2]{I2021} can be proven as in \cite{I2021} invoking Theorem~\ref{general-theorem} in place of \cite[Theorem~1.3]{I2021}.  However, we can actually obtain a stronger result, which we state here, in that we do not need the hypothesis made in \cite[Theorem~1.2]{I2021} that $\Omega$ is a Runge domain of holomorphy.  (Note the \cite[Theorem~1.3]{I2021} involved a Runge domain of holomorphy but that this is not the case with Theorem~\ref{general-theorem} above.)

\bthm\label{poly-convex}
Let $\Omega$ be a bounded, connected open set in $\C^N$, 
let $x_0$ be a point of $\Omega$, and let $\vep>0$.  
Then there exists a polynomially convex arc $J$ in $\C^N$ such that $x_0\in J \subset \Omega$ and $m(\Omega\sm J)<\vep$.  Furthermore, $J$ can be chosen so that $P(J)=C(J)$ and the set of polynomials zero-free on $J$ is dense in 
$P(\overline\Omega)$.  The same statements hold with \lq\lq arc\rq\rq\ replaced by \lq\lq simple closed curve\rq\rq\ provided $N\geq 2$.  
\ethm

The proof of Theorem~\ref{poly-convex} is the same as that given for \cite[Theorem~1.2]{I2021} except for invoking Theorem~\ref{general-theorem} above in place of \cite[Theorem~1.3]{I2021} and replacing \cite[Theorem~2.4]{I2021} by the more general Lemma~\ref{big-set} below.

In the next section we present some lemmas.  Theorem~\ref{general-theorem} and Corollary~\ref{corollary} are proved in Section~3.

%
%
%
%

\section{Lemmas}

We begin with some topological lemmas concerning arcs.

\blem\label{prelemma}
Let $\lambda$ be a closed set in $\R^N$, $N\geq 3$, of topological dimension at most $1$, let $a$ and $a'$ are two points in $\lambda$, and let $\Omega$ be a connected open set of $\R^N$ that contains $a$ and $a'$.  Then there is an arc $J$ from $a$ to $a'$ contained in $\Omega$ that intersects $\lambda$ only in the end points $a$ and $a'$ of $J$ and is such that the open arc $J\sm \{a,a'\}$ is $C^\infty$-smooth.
\elem

\bpf
The proof is essentially a repetition of the 
proof of the $n\geq 3$ case of \cite[Theorem~1.2]{IS}.
\epf

The above lemma becomes false with $N=2$.  There is, however, the following weaker result.

\blem\label{prelemma-2}
Let $\lambda$ be an arc in $\R^2$, let $a$ and $a'$ are two points in $\lambda$, let $\lambda_{a,a'}$ be the subarc of $\lambda$ from $a$ to $a'$, and let $\Omega$ be a neighborhood of $\lambda_{a,a'}$ in $\R^2$.  Then there is an arc $J$ from $a$ to $a'$ contained in $\Omega$ that intersects $\lambda$ only in the end points $a$ and $a'$ of $J$ and is such that the open arc $J\sm \{a,a'\}$ is $C^\infty$-smooth.
\elem

\bpf
The proof is similar to the 
proof of the $n=2$ case of \cite[Theorem~1.2]{IS} but requires some care, so we include the details.  There is a conformal map $\varphi:\U\rightarrow\C^*\sm\lambda$
of the open unit disc $\U$ onto the complement of $\lambda$ in the Riemann sphere $\C^*$.  By \cite[Theorem~2.1]{Church} (which seems to be due to Marie Torhorst \cite{Torhorst}), $\varphi$ extends continuously to the closure of $\U$.  Choose points $p$ and $p'$ in the boundary $\partial \U$ of $\U$ such that $\varphi(p)=a$ and $\varphi(p')=a'$ and such that on one of the open arcs $\alpha$ on $\partial \U$ determined by $p$ and $p'$ the function $\varphi$ never takes either of the values $a$ and $a'$.  Then $\varphi$ maps $\alpha$ onto the interior of the arc $\lambda_{a,a'}$.  Let $\ell$ be an arc in $\U\cup \{p,p'\}$ with end points $p$ and $p'$ and interior contained in $\U$, and set $J=\varphi(\ell)$.
Then $J$ is an arc from $a$ to $a'$ that intersects $\lambda$ only in the points $a$ and $a'$.
By choosing $\ell$ to lie sufficiently near $\alpha$ and to be $\C^\infty$-smooth, we can insure that $J$ is contained in $\Omega$ and that $J\sm\{a,a'\}$ is 
$C^\infty$-smooth.
\epf

\blem\label{connected}
Let $\Omega$ be a connected open set in $\R^2$, and let $J_1,\ldots,J_n$
be finitely many disjoint arcs in $\Omega$.  Then $\Omega\sm (J_1\cup\cdots\cup J_n)$ is connected.
\elem

\bpf
We merely sketch the proof leaving the details to the reader.  By induction it suffice to consider the case when there is just one arc $J$.
For that case, first show that $\Omega$ contains a connected neighborhood $V$ of $J$ such that the complement of $V$ in $\R^2$ is connected.  Since $V$ is then homeomorphic to the plane, it is a standard fact that $V\setminus J$ is connected.  Connectedness of $\Omega\setminus J$ follows.
\epf

\blem\label{smooth-arc}
Let $E$ be a compact set that is contained in an arc in $\R^N$.  Then every connected neighborhood of $E$ contains an arc $J$ that contains $E$ and has the additional property that each component of $J\sm E$ is a $C^\infty$-smooth open arc.
\elem

\bpf
The case $N=1$ is trivial.  We first treat the case $N\geq 3$; the case $N=2$ requires a more involved argument.

Let $\Omega$ be a connected neighborhood of $E$.  
Let $\sigma: \I\rightarrow \R^N$ be an arc through $E$, and assume without loss of generality that the end points of $\sigma$ are in $E$.  Set $L=\sigma^{-1}(E)$.  
The set $\I\sm L$ is an at most countable union of disjoint open intervals $\openab 1, \openab 2, \ldots$.  Note that the diameters $\diam\bigl(\sigma(\ab j)\bigr)$ go to zero as $j\rightarrow\infty$ (if there are infinitely many intervals $\openab j$).  In particular, $\sigma(\ab j)$ is contained in $\Omega$ for all but finitely many $j$.  For each $j$ such that $\sigma(\ab j)$ is contained in $\Omega$, choose a connected neighborhood $\Omega_j$ of $\sigma(\ab j)$ contained in $\Omega$ and of diameter no more than twice the diameter of $\sigma(\ab j)$.  For $j$ such that $\sigma(\ab j)$ is not contained in $\Omega$, set $\Omega_j=\Omega$.  Note that then $\diam(\Omega_j)\rightarrow 0$ and $j\rightarrow \infty$.

By Lemma~\ref{prelemma} there is an arc $\sigma_1:\ab 1\rightarrow \Omega_1$ from $\sigma(a_1)$ to $\sigma(b_1)$ that is $C^\infty$-smooth except possibly at its end points and that intersects $E$ only in its end points.  Continuing inductively, we can choose, for each $j=2,3,\ldots$, an arc $\sigma_j:\ab j\rightarrow \Omega_j$ from $\sigma(a_j)$ to $\sigma(b_j)$ that is $C^\infty$-smooth except possibly at its end points and that intersects 
$E\cup \sigma_1(\ab 1)\cup \cdots\cup \sigma_{j-1}(\ab {j-1})$ only in its end points.  Now defining $\tau:\I\rightarrow\Omega$ to coincide with 
$\sigma$ on $L$ and to coincide with $\sigma_j$ on $\ab j$ for each $j=1,2,\ldots$ yields the desired arc.  (Continuity of $\tau$ is a consequence of the conditions that $\diam\bigl(\sigma_{j}(\ab j)\bigr)\leq 2\,\diam(\Omega_j)$ and $\diam(\Omega_j)\rightarrow 0$ as $j\rightarrow\infty$.)
This concludes the proof in the case $N\geq 3$. 

We now consider the case $N=2$, which we will establish in two steps.  First we will obtain an arc through $E$ that is contained in $\Omega$ and for which no smoothness is asserted, and subsequently we will obtain the arc whose existence is asserted in the statement of the lemma.

Let $\Omega$, $\sigma$, and $L$ be as before. 
Since $L$ is a compact set contained in the (relatively) open set $\sigma^{-1}(\Omega)$ of $[0,1]$, the set $L$ is contained in a finite union of intervals open in $[0,1]$ and contained in $\sigma^{-1}(\Omega)$.
Consequently, we can choose points 
$$0=c_0<d_0<c_1<d_1<\cdots<c_n<d_n=1$$ such that 
$$E\subset \sigma\bigl([c_0,d_0] \cup\cdots\cup [c_n,d_n]\bigr) 
\subset \Omega.$$
By modifying $\sigma$ near each of the points $d_0,\ldots,d_{n-1}$ and $c_1,\ldots, c_n$, we may assume that there are open Euclidean balls $B_{d_0},\ldots,B_{d_{n-1}}$ and $B_{c_1},\ldots, B_{c_n}$  centered at $\sigma(d_0),\ldots, \sigma(d_{n-1})$ and
$\sigma(c_1),\ldots, \sigma(c_n)$, respectively, whose closures are disjoint and lie in $\Omega$ such that the intersection of $\sigma(\I)$ with each of these balls is a straight line segment.  Choose points $q_0',\ldots, q_{n-1}'$ and $p_1',\ldots, p_n'$ in $B_{d_0},\ldots,B_{d_{n-1}}$ and $B_{c_1},\ldots, B_{c_n}$, respectively.  
The set 
$\Omega\sm \sigma\bigl([c_0,d_0] \cup\cdots\cup [c_n,d_n]\bigr)$ is connected by Lemma~\ref{connected}, so there is an arc from $q_0'$ to $p_1'$ in $\Omega\sm \sigma\bigl([c_0,d_0] \cup\cdots\cup [c_n,d_n]\bigr)$.
By discarding initial and final segments of this arc, we can obtain an arc $\lambda$ in
$\Omega\sm \bigl[\sigma\bigl([c_0,d_0] \cup\cdots\cup [c_n,d_n]\bigr) \cup B_{d_0} \cup B_{c_1}\bigr]$
whose end points $\tilde q$ and $\tilde p$ lie on the boundary of $B_{d_0}$ and $B_{c_1}$, respectively.
Let $L_{\tilde q}$ and $L_{\tilde p}$ be the straight line segments from 
$\sigma(d_0)$ to $\tilde q$ and from $\sigma(c_1)$ to $\tilde p$, respectively.  Set $\ell_0=L_{\tilde q} \cup \lambda \cup L_{\tilde p}$. Then $\ell_0$ is an arc in 
$\Omega$ that intersects $\sigma\bigl([c_0,d_0]  \cup\cdots\cup [c_n,d_n]\bigr)$ only in the end points $\sigma(d_0)$ and $\sigma(c_1)$ of $\ell_0$.

Continuing inductively we can similarly choose, for each $j=1,\ldots, n-1$, an arc $\ell_j$ in $\Omega$ from $\sigma(d_j)$ to $\sigma(c_{j+1})$ that, aside from its end points, is disjoint from 
$\sigma\bigl([c_0,d_0] \cup\cdots\cup [c_n,d_n]\bigr)$ and from each of the previously chosen arcs $\ell_1,\ldots,\ell_{j-1}$.  Then 
$\sigma\bigl([c_0,d_0] \cup\cdots\cup [c_n,d_n]\bigr) \cup \ell_0\cup\cdots\cup\ell_{n-1}$ is an arc in $\Omega$ that contains $E$.

The passage from the arc just obtained to the desired arc $J$ is similar to the proof of the lemma in the case $N\geq 3$ but with Lemma~\ref{prelemma} replaced by Lemma~\ref{prelemma-2}, so we compress the details.  Let $\tsigma:[0,1]\rightarrow\Omega$ be a parametrization of the arc just obtained satisfying $\tsigma|L=\sigma|L$.  With $[0,1]\sm L = (a_1, b_1) \cup (a_2, b_2) \cup \cdots$ as in the case $N\geq3$, we choose, for each $j=1,2,\ldots$, a connected neighborhood $\Omega_j$ of $\sigma(\ab j)$ contained in $\Omega$ in such a way that $\diam(\Omega_j)\rightarrow 0$ as $j\rightarrow\infty$.  By Lemma~\ref{prelemma-2} there is an arc $\sigma_1:\ab 1\rightarrow \Omega_1$ from $\tsigma(a_1)$ to $\tsigma(b_1)$ that is 
$C^\infty$-smooth except possibly at its end points and that intersects $\tsigma([0,a_1]\cup[b_1,1])$ only in its end points $\tsigma(a_1)$ and $\tsigma(b_1)$.
Define $\tau_1:[0,1]\rightarrow \Omega$ to coincide with $\tsigma$ on $[0,a_1]\cup[b_1,1]$ and to coincide with $\sigma_1$ on $\ab 1$.  Then $\tau_1$ is an arc.
Continuing inductively, we can obtain, for each $j=2,3,\ldots$, an arc $\sigma_j:\ab j\rightarrow \Omega_j$ from $\tsigma(a_j)$ to $\tsigma(b_j)$ and an arc $\tau_j$ so that $\sigma_j$ is $C^\infty$-smooth except possibly at its end points and intersects $\tau_{j-1}([0,a_j]\cup[b_j,1])$ only in its end points $\tsigma(a_j)$ and $\tsigma(b_j)$, and $\tau_j$ coincides  with $\tau_{j-1}$ on $[0,a_j]\cup[b_j,1]$ and coincides with $\sigma_j$ on $\ab j$.
The sequence $(\tau_n)$ converges uniformly, and its limit is the desired arc.
\epf

Our final two lemmas concern polynomial convexity.  The first of these
is an almost immediate consequence of \cite[Lemma~3.2]{Iaccepted}, and as mentioned in the introduction generalizes \cite[Theorem~2.4]{I2021}.

\blem\label{big-set}
Let $Y$ be a compact set in $\C^N$, let $x_0$ be a point of $Y$, and let $\vep>0$.
Let $\{p_j\}$ be a countable collection of polynomials on $\C^N$ such that $p_j(x_0)\neq 0$ for all $j$.  Then there exists a totally disconnected, compact polynomially convex set $E$ with $x_0\in E\subset Y\subset \C^N$ such that
\item{\rm(i)} each $p_j$ is zero-free on $E$
\item{\rm(ii)} $m(Y\sm E)<\vep$.  
\elem

\bpf 
By making a complex affine change of coordinates, we may assume without loss of generality that $x_0=0$ and that $Y$ is contained in the open unit ball $B$ of $\C^N$.  Then \cite[Lemma~3.2]{Iaccepted} gives a totally disconnected, compact polynomially convex set $K$ with $0\in K\subset B\subset \C^N$ such that
\smallskip
\item{\rm(i)} each $p_j$ is zero-free on $K$
\item{\rm(ii)} $m(B\sm K)<\vep$.
\smallskip\hfil\break
Let $E=K\cap Y$. Then $E$ is compact and totally disconnected.  Also $x_0=0\in E\subset Y\subset\C^N$, each $p_j$ is zero-free on $E$, and $m(Y\sm E)\leq m(B\sm K)<\vep$.  Because $K$ is a totally disconnected, compact polynomially convex set, it follows from the Shilov idempotent theorem that $P(K)=C(K)$ (see \cite[p.~48, Corollary~3]{AW} for instance), and hence, $P(E)=C(E)$, so $E$ is polynomially convex.
\epf

The following lemma, a special case of \cite[Corollary~1.4]{IS}, will play a key role in the proof of Theorem~\ref{general-theorem}.

\blem\label{Stout-lemma}
Let $Y\subset\C^N$ be a compact polynomially convex set, and let $\Gamma$ be a rectifiable arc both of whose end points lie in $Y$ but is otherwise disjoint from $Y$.  If a nonempty open subarc of $\Gamma$ is contained in a purely one-dimensional analytic subvariety $V$ of $\C^N$ but $\Gamma$ is not entirely contained in $V$, then $Y\cup\Gamma$ is polynomially convex.
\elem

%
%
%
%

\section{Proofs of Theorem~\ref{general-theorem} and Corollary~\ref{corollary}}

Corollary~\ref{corollary} is a special case of Theorem~\ref{general-theorem}, but we give an independent proof since the corollary is much more easily established than the general theorem.

\begin{proof}[Proof of Corollary~\ref{corollary}]
Let $\Omega$ be an arbitrary connected neighborhood of $E$.
Let $J$ be the arc in $\Omega$ given by Lemma~\ref{smooth-arc}.  Then \cite[Theorem~1.7]{IS} gives that $J$ is polynomially convex and that $P(J)=C(J)$ if $P(E)=C(E)$.  The statement about a simple closed curve then follows from \cite[Theorem~1.1]{IS}.
\end{proof}

\begin{proof}[Proof of Theorem~\ref{general-theorem}]
The case $N=1$ follows immediately from Lemma~\ref{smooth-arc} since every arc is polynomially convex.  From now on we assume that $N\geq 2$. 
We treat only the construction of the arc. The simple closed curve can be constructed similarly, or alternatively, it can be obtained from the arc by invoking \cite[Theorem~1.1]{IS}. 

Let $\Omega$ be an arbitrary connected neighborhood of $E$.
Lemma~\ref{smooth-arc} yields the existence of an arc $\sigma_0:\I\rightarrow \Omega$ such that $\sigma_0(\{0,1\})\subset E\subset\sigma_0(\I)$ and such that the restriction $\sigma_0|_{(\I\sm\sigma_0^{-1}(E))}$ is a $C^\infty$-embedding.  The proof will be complete once we establish that there is an arc $\sigma:\I\rightarrow\Omega$ that has the properties just listed for $\sigma_0$ and has the additional property that setting $J=\sigma(\I)$ we have that $\h J = J \cup \h E$.  We will obtain $\sigma$ as a uniform limit of a sequence of arcs $\sigma_n:\I\rightarrow \Omega$ that will be constructed inductively.  We will also simultaneously construct a sequence $(K_n)$ of compact polynomially convex sets.

Let $K_0$ be a closed ball in $\C^N$ whose interior contains $\sigma_0(\I)$.
Choose a sequence $(U_n)$ of neighborhoods of $\h E$ with $U_0=\C^N$ and $\bigcap_{n=0}^\infty U_n= \h E$. 
Set $L=\sigma_0^{-1}(E)$.
Roughly, the $\sigma_n$ will be constructed by succesively changing $\sigma$ at most once on each component of $\I\sm L$.

The sequence of arcs $(\sigma_n)$ and the sequence of compact polynomially convex sets $(K_n)$ will be chosen so that the following conditions hold for all $n=1,2,\ldots$.
\itemskip
{
\begin{enumerate}
\item[(i)] $\h E\subset \Int(K_n)\subset K_n\subset \Int(K_{n-1})\cap U_n$.
\itemskip
\item[(ii)] $K_n$ has $C^\infty$-smooth boundary $\partial K_n$.
\itemskip
\item[(iii)]  $\sigma_n|_L=\sigma_0|_L$.
\itemskip
\item[(iv)]  $\| \sigma_n - \sigma_{n-1}\|_\infty < 1/n$.  (Here $\|\cdot\|_\infty$ denotes the supremum~norm.)
\itemskip
\item[(v)] The restriction of $\sigma_n$ to $\I\sm L$ a $C^\infty$-immersion.
\itemskip
\item[(vi)] $\sigma_{n-1}$ is transverse to $\partial K_n$.
\itemskip
\item[(vii)] $\sigma_n(\I)\subset K_{n-1}\cup \sigma_{n-1}(\I)$. 
\itemskip
\item[(viii)] On account of condition (vi), $\sigma_{n-1}^{-1}(\C^N\sm K_n)$ is a finite union of open intervals $\openab 1,\ldots, \openab u$.  The map $\sigma_n$ coincides with $\sigma_{n-1}$ everywhere on $\I$ expect on those intervals $(a_j,b_j)$ such that $\sigma_{n-1}(a_j,b_j)$ is contained entirely in $K_{n-1}$.  Also $\sigma_n\bigl(\openab 1\cup\cdots\cup\openab u\bigr)$ is disjoint from $K_n$.
\itemskip
\item[(ix)]  For each component $\lambda$ of $\sigma_n(\I)\sm K_n$ there is an analytic subvariety of $\C^N$ that contains a nonempty open subarc of $\lambda$ but does not contain all of~$\lambda$.
\end{enumerate}
}

Before constructing the sequences $(\sigma_n)$ and $(K_n)$ we show that their existence will yield the theorem.  The set $\I\sm L$ consists of an at most countable collection of open intervals $(\alpha_1,\beta_1),(\alpha_2,\beta_2),\ldots$.  From condition (viii) we get that the sequence $(\sigma_n)$ has a very special form.  If $\sigma_0\bigl(\albeta\bigr)$ is contained in $\h E=\bigcap_{n=0}^\infty U_n=\bigcap_{n=0}^\infty K_n$, then every $\sigma_n$ coincides with $\sigma_0$ on $\albeta$.  Otherwise there is a smallest integer $s\geq 1$ such that $\sigma_s\bigl(\albeta\bigr)$ intersects the complement of $K_s$, and then $\sigma_{s-1}$ and $\sigma_s$ may differ on $\albeta$, but $\sigma_0,\ldots,\sigma_{s-1}$ coincide on $\albeta$, and $\sigma_n$ coincides with $\sigma_s$ on $\albeta$ for all $n\geq s$.  

In combination with conditions (iii) and (iv), the observation in the preceding paragraph yields that the sequence $(\sigma_n)$ converges uniformly to a map $\sigma:\I\rightarrow\Omega$.  Furthermore, $\sigma$ is injective because given $t_1\neq t_2$ in $\I$ there exists an $s$ such that $\sigma(t_1)=\sigma_s(t_1)\neq \sigma_s(t_2)=\sigma(t_2)$.  Thus $\sigma$, or more precisely $J=\sigma(\I)$, is an arc.  Also since for each $j=1,2,\ldots$ there is an $m$ such that $\sigma$ coincides with $\sigma_m$ on $\albeta$, condition (v) yields that each component of $J\sm E$ is a $C^\infty$-smooth open arc.

Condition (iii) implies that $\sigma|L=\sigma_0|L$, so $J\supset \sigma_0(L)=E$.  Consequently, $\h J\supset J\cup \h E$.  To establish the reverse inclusion it suffices to show that $J\cup \h E$ is polynomially convex.

By conditions (viii) and (ix), there are finitely many disjoint arcs $\lambda_1,\ldots,\lambda_u$ such that $K_n\cup \sigma_n(\I)= K_n\cup \lambda_1\cup\cdots\cup\lambda_u$ where each $\lambda_j$ intersects $K_n$ precisely in its two end points and is such that there is an analytic subvariety of $\C^N$ that contains a nonempty open subarc of $\lambda_j$ but does not contain all of $\lambda_j$.  Consequently, successive application of Lemma~\ref{Stout-lemma} shows that $K_n\cup \sigma_n(\I)$ is polynomially convex.  We will show that 
\begin{equation}\label{intersection-eq}
\textstyle\bigcap\limits_{n=0}^\infty\Bigl[K_n\cup \sigma_n(\I)\Bigr]=J\cup\h E
\end{equation}
thereby establishing the polynomial convexity of $J\cup\h E$.
Clearly $\h E$ is contained in the left hand side of equation~(\ref{intersection-eq}) by condition (i).  
For each point $x_0$ in $J$ there is an $s$ such that $x_0$ is in $\sigma_n(\I)$ for all $n\geq s$.  Since conditions conditions (i) and (vii) show that the sequence of sets $\bigl(K_n\cup\sigma_n(\I)\bigr)$ is decreasing, this yields that $J$ is contained in the left hand side of equation~(\ref{intersection-eq}) also.
Thus the left hand side of equation~(\ref{intersection-eq}) contains the right hand side.
For the reverse inclusion note that for a point $x_0$ in the left hand side of equation~(\ref{intersection-eq}) that does not lie in $\h E$ there is an $s$ such that $x_0$ does not lie in $K_s$; then $x_0$ is in $\sigma_n(\I)$ for all $n\geq s$, and thus $x_0$ must be in $J$.  This concludes the verification that the existence of the sequences $(\sigma_n)$ and $(K_n)$ will yield the theorem.

It remains to construct the sequences $(\sigma_n)$ and $(K_n)$.  We already have $\sigma_0$ and $K_0$.
We proceed now by induction.  Suppose for some $k\geq 0$ we have chosen
$\sigma_0,\ldots, \sigma_k$ and $K_0,\ldots, K_k$ so that conditions (i)--(ix) are satisfied for all $n=1,\ldots, k$.  By well-known theorems in several complex variables regarding the existence of plurisubharmonic exhaustion functions and the equality of polynomial hulls and plurisubharmonic hulls (see for instance \cite[Theorems~II.5.11 and~VI.1.18]{Range}) there exists a $C^\infty$ strictly plurisubharmonic exhaustion function $\varphi$ on $\C^N$ such that
\vskip-5pt
\[
\varphi(z)>0\ {\rm for\ } z\in \C^N\sm (\Int(K_k)\cap U_ {k+1})
\]
\vskip-8pt
and
\[
\varphi(z)<0\ {\rm for\ } z\in E.\hskip 100pt
\]
\vskip8pt
Set $K^r=\{z\in \C^N: \varphi(z)\leq r\}$.  Set $M=\max_{z\in E} \varphi(z)$.  Then for $M<r<0$ the set $K^r$ is a compact polynomially convex set such that
$\h E\subset \Int(K^r)\subset K^r\subset \Int(K_k)\cap U_{k+1}$.  Furthermore, since for these values of $r$ the boundary of $K^r$ is disjoint from $E$, we can choose, by Sard's theorem, a value $r_0$ satisfying $M<r_0<0$ such that $K^{r_0}$ has $C^\infty$-smooth boundary $\partial K^{r_0}$ and $\sigma_k$ is transverse to $\partial K^{r_0}$.  Set $K_{k+1}=K^{r_0}$.

Because $\sigma_k$ is transverse to $K_{k+1}$, the set $\sigma_k^{-1}(\C^N\sm K_{k+1})$ is a finite union of open intervals.  Of those open intervals let $\openab 1,\ldots,\openab w$ denote those whose image under $\sigma_k$ is entirely contained in $K_k$.  For each $j=1,\ldots, w$, choose a point $x_j$ in $\sigma_k\bigl(\openab j\bigr)\cap \Int(K_k)$.  Then choose disjoint closed balls $\ob_1,\ldots,\ob_w$ centered at $x_1,\ldots, x_w$, respectively, with radii strictly less than $1/2(k+1)$ and small enough that each $\ob_j$ is contained in $K_k$ and the intersection of each $\ob_j$ with $K_{k+1} \cup \sigma_k(\I)$ is contained in $\sigma_k\bigl(\openab j\bigr)$.  Then choose, for each $j=1,\ldots, w$, a $C^\infty$-smooth arc $\tau_j:\ab j\rightarrow \Omega$ such that 
$\tau_j$ coincides with $\sigma_k|_{\ab j}$ except on some subarc of $\openab j$ that is mapped into $\ob_j$ by both $\sigma_k$ and $\tau_j$, and such that there is an analytic subvariety of $\C^N$ that contains an open subarc of $\tau_j(\ab j)$ but does not contain all of $\tau_j(\ab j)$.  Finally, define $\sigma_{k+1}:\I\rightarrow \Omega$ to coincide with $\sigma_k$ on $\I\sm \bigl(\openab 1\cup\cdots\cup\openab w\bigr)$ and to coincide with $\tau_j$ on $\ab j$ for each $j=1,\ldots, w$.  Then $\sigma_{k+1}$ is an arc.  Furthermore, conditions (i)--(ix) hold for all $n=1,\ldots, k+1$.  This completes the induction and the proof.
\end{proof}

\end{document}